\documentclass[11pt,a4]{article}
\usepackage{amssymb}
\usepackage{amsfonts}
\usepackage{amsthm}
\usepackage{amsbsy}
\usepackage{graphicx}
\usepackage{verbatim}
\usepackage{enumerate}
\usepackage{enumitem}
\usepackage{hyperref}
\def\BE{\begin{equation}}
\def\EE{\end{equation}}
\def\BEA{\begin{eqnarray}}
\def\EEA{\end{eqnarray}}

\def\NNL{\nonumber \\}
\def\NN{\nonumber}
\def\NewAND{\hspace{-0.15in} & & \hspace{-0.15in}}

\newtheorem{Theo}{Theorem}[section]
\newtheorem{Lem}{Lemma}[section]
\newtheorem{Def}{Definition}[section]

\newtheorem{Rem}{Remark}[section]

\theoremstyle{Definition} 
\providecommand{\bysame}{\leavevmode\hbox to3em{\hrulefill}\thinspace}
\providecommand{\MR}{\relax\ifhmode\unskip\space\fi MR }
\begin{document}

\title{Hidden Markov Mixture Autoregressive Models:
Stability and Moments}

\author{S.H.Alizadeh, S.Rezakhah\footnotemark[1]}
\footnotetext[1]{	
            Faculty of Mathematics and Computer Science
			,Amirkabir University of Technology, Tehran, Iran.
            Email:\texttt{rezakhah@aut.ac.ir,\, sasan\_alizadeh@aut.ac.ir} \par
			}
\date{}
\maketitle

\begin{abstract}
\textbf{T}his paper \textbf{introduces} a new parsimonious structure for mixture of autoregressive model\textbf{s}. \textbf{The} weighting coefficients are determined through latent random variables, follow\textbf{ing} a hidden Markov model. \textbf{We propose a} dynamic programming algorithm for the application of forecasting. We also derive the limiting behavior of  unconditional first moment of the process and an appropriate upper bound for the limiting value of the variance. \textbf{This} can be considered as long run behavior of the  process. Finally we show convergence and stability of the second moment. \textbf{Further, we illustrate the efficacy of the proposed model by simulation and forecasting.}
\end{abstract}


 \par

\section{Introduction}

The most frequently used approaches to time series model building assume
that the data under study are generated from a linear stochastic process.
Linear models provide
a number of appealing properties (such as physical interpretations, frequency
domain analysis, asymptotic results, statistical inference and many others)\cite{Brockwell}.
Despite those advantages, it is well known that real-life systems are usually nonlinear, and certain
features, such as limit-cycles, asymmetry \cite{Lai},\cite{Medeiros}, conditional heteroscedasticity \cite{Engle}, flat
stretches, bursts  \cite{LMR} and jump phenomena cannot be correctly captured by linear statistical models.

Since the Mixture Transition Distribution (MTD) was originally introduced by Raftery \cite{Raftery} for modeling high order Markov chains in the discrete state space, the broad family of this model have been extended and applied for modeling conditional distribution of observations in the context of nonlinear time series with arbitrary state spaces \cite{BerchRaft}. This model also has been extended to the mixture transition of Gaussian distributions, known as GMTD, which contains autoregressive model as a special case, for modeling flat stretches, bursts and outliers \cite{LMR}.
Mixture of Autoregressive (MAR) model \textbf{(}which has been proposed by Wong and Li \cite{WL}\textbf{)} is a flexible generalization of GMTD to model processes with multimodal conditional distributions and conditional heteroscedasticity.
The important feature of MAR model is that it can be considered as the mixture of some stationary  and non-stationary AR processes and remains stationary.
For time series $\{Y_t\}_{t=0}^\infty$,  $Y_t \in \mathbb{R}$, the MAR($K;p_1,p_2,...,p_K$) is defined as
\BEA
    F(y_t|\mathcal{F}_{t-1}) = \sum_{k=1}^K \alpha_k \Phi(\frac{y_t-\phi_{k0}-\phi_{k1}y_{t-1}-...-\phi_{k{p_k}}y_{t-{p_k}}}{\sigma_k}), \label{MAR}
\EEA
in which  $y_t$ denotes a realization of $Y_t$, and $\mathcal{F}_{t}=\sigma\{Y_s:s \leq t\}$ and $F(y_t|\mathcal{F}_{t-1})$ is the \textbf{conditional} distribution of $Y_t$ given information of $\mathcal{F}_{t-1}$. Also  $\alpha_k,\,\, k=1,\cdots,K$ are the weighting coefficients (i.e. $\alpha_k>0,\,k=1,...,K$ and $\sum_{k=1}^K \alpha_k = 1$.)  and $\Phi(.)$ is the cumulative distribution function of the standard normal distribution. This model  is a mixture of $K$ Gaussian AR($p_k$), $k=1,...,K$ models \cite{WL}.

The mixture of autoregressive conditional heteroscedasticity model was also proposed by Wong and Li \cite{WL2001} to capture the squared autocorrelation structure of observations.
Berchtold \cite{Berchtold} also introduced a new approach for modeling heteroscedastic time series with MTD model in which the variances of each Gaussian distributions depends on the past time series observations. For exhaustive review of MTD model see \cite{BerchRaft}.

In  the MTD models the contribution of distributions are always fixed and it is not sensitive to the past observations. However for real processes one might expect better forecast interval if additional information from the past were allowed to affect \cite{Engle}.
Another approach to study mixture models is to introduce some latent variables $\{Z_t\}_{t=p+1}^\infty$, which are iid and $Y_t$ given $Z_t$ is independent of $\{Z_s:\,s \neq t\}$.
Each variable $Z_t$ has a discrete distribution with support $\{1,\cdots,K\}$ \textbf{with} probability \textbf{masses} $P(Z_t=k)=\alpha_k,\, k=1,\cdots,K$ as the weighting coefficients in the mixture model.
Since these models \textbf{do} not consider the dependency structure of latent variables, the dynamics of weighting coefficients can not be modeled.
For finite state space time series, Bartolucci and Farcomeni \cite{Bartolucci} studied a generalization of mixture transition models with hidden Markov models.

In this paper\textbf{,} we propose a new approach to model conditional distribution of $Y_t$ given past information for nonlinear time series in general state space (i.e.  $Y_t \in R$). We use latent Markov process as an appropriate tool to consider the effect of past information and build a parsimonious model; the idea of Markov switching models (see Hamilton \cite{Hamilton}, Mcculloch and Tsay \cite{Tsay}) for process $\{Y_t\}_{t=0}^\infty$.
Our new model include\textbf{s} the hidden Markov model (HMM) \cite{Bishop} as a special case and it also generalizes MAR model in a sensible way.
This model makes use of the whole past information to maximize the posterior probability of $Z_{t-1}$ \textbf{(}given observed $Y_0,\cdots,Y_{t-1}$\textbf{)} and predicts the probability of  $Z_{t}$ by the Markov assumption of the latent process. Although using \textbf{all} past \textbf{observations could increase the complexity of} the model, we propose a dynamic programming algorithm which reduces the volume of calculations for  \textbf{forecasting}.
We derive the limiting behavior of the first unconditional moment of the process, and obtain an upper bound for the limit of variance. \textbf{We also} investigate the existence and stability of the second moment.

\textbf{This paper is organized as follows. Hidden Markov Mixture Autoregressive (HM-MAR) model is introduced in section \ref{Sec HM-MAR}. Section \ref{Sec StatProp} is devoted to the statistical properties of the HM-MAR model. Section \ref{Sec Simulation} analyzes the efficiency of the proposed model through simulation and comparison of the forecast errors with the MAR model. Section \ref{Sec Conclusion} concludes the paper.}

\section{Hidden Markov Mixture Autoregressive Model} \label{Sec HM-MAR}

Let $Y=\{Y_t\}_{t=0}^\infty$ be a sequence of random variables in $\mathbb{R}$ \textbf{where}  $y_t$ \textbf{is} a realization of $Y_t$. \textbf{Also let} $\mathcal{F}_{t}=\sigma\{Y_s:s \leq t\}$ \textbf{and $F(y_t|\mathcal{F}_{t-1})$ respectively represent} the sigma-field of all information up to time $t$, and  the conditional distribution function of $Y_t$ \textbf{(}given past information and  $\alpha_h^{(t)} \equiv \alpha_h^{(t)}(y_1,...,y_{t-1})$\textbf{)}. \textbf{In addition}
 $\{Z_t\}_{t \geq p}$ \textbf{denotes} a hidden or latent process,  \textbf{a} positive recurrent Markov chain on a finite set $E=\{1,2,...,K\}$\textbf{.} \textbf{The} initial conditional probabilities \textbf{are}
\BEA
   \boldsymbol{\rho}=(\rho_1,\cdots,\rho_K)^\prime,\,\,\, \rho_h=P(Z_{p}=h|y_0,\cdots,y_{p-1}) \qquad
   h=1,...,K, \label{Initial Probabilities}
\EEA
\textbf{with} transition probability matrix
\BEA
    P=\|\pi_{i,j}\|_{K \times K},
    \label{Transition Probability Matrix}
\EEA
in which
\BEA
    \pi_{i,j}=P(Z_t=j|Z_{t-1}=i),\qquad i,j \in \{1,...,K\},
    \label{Transition Probabilities}
\EEA
and invariant probability measure
\BEA
    \boldsymbol{\mu}=(\alpha_1,...,\alpha_K)^\prime, \label{Mu}
\EEA
\textbf{where} $\alpha_j=\lim_{t \rightarrow \infty} P(Z_t=j)$.

We consider $\{Y_t\}_{t=0}^\infty$ to have a Hidden Markov-Mixture Autoregressive, HM-MAR($K,p$), model with $K$ normal distributions, $p$ lagged observations in the AR processes, if the conditional distribution of $Y_t$ given $\mathcal{F}_{t-1}$ is defined as follows:

\begin{enumerate}
    \item[i.]

        \textbf{For} $t=p$
         \BEA
            F(y_{p},Z_p=h|\mathcal{F}_{p-1}) =  \rho_h
                                   \Phi(\frac{y_p-a_{0,h}-a_{1,h}y_{p-1}-...-a_{{p},h}y_{0}}{\sigma_h}),
                                   \label{HM MAR p1}
         \EEA

    \item[ii.]

        \textbf{For} $t \geq p + 1$
        \BEA
            F(y_{t}|\mathcal{F}_{t-1})= \sum_{h=1}^K \alpha_h^{(t)}
            \Phi(\frac{y_t-a_{0,h}-a_{1,h}y_{t-1}-...-a_{{p},h}y_{t-{p}}}{\sigma_h}),
            \label{HM MAR}
        \EEA
\end{enumerate}
\textbf{where} $\alpha_h^{(t)}=P(Z_t=h|\mathcal{F}_{t-1})$ and $\Phi(.)$ is the standard normal distribution function.

In fact latent random variables $\{Z_t\}_{t=p+1}^\infty$ determine the contribution of distributions in the mixture model and conditioning on $Z_t$\textbf{.} \textbf{We} assume  $Y_t$ is $p$-tuple Markov\textbf{,} independent of $\{Z_{s},\, s \neq t \}$. In other words\textbf{, by} conditioning on $\{Y_{t-1},\cdots,Y_{t-p}\}$ and $Z_t$, $Y_t$ \textbf{is independent of} $\{Y_{s}, \, s < t-p\}$ and $\{Z_{s},\, s \neq t \}$.

The novelty of HM-MAR model is that
the contribution of each distribution in the mixture structure is not of predefined fixed form.
It makes use of the \textbf{all} past observations from $Y_0$ up to $Y_{t-1}$.
The hidden Markov assumption of the process $\{Z_t\}_{t \geq p}$, enables us to build a parsimonious model.

The MAR model \cite{WL} can be considered as a special case of such a HM-MAR model (\ref{HM MAR p1}-\ref{HM MAR}), in which the transition matrix $P$ of the process $\{Z_t\}_{t \geq p}$ has $K$ identical rows (i.e. $p(Z_t=i|Z_{t-1}=j)=\alpha_i$ for all $i,\,j=1,...,K$\textbf{. T}hat is $\{Z_t\}_{t=p+1}^\infty$ are independent and identically distributed) with $p(Z_t=i|Z_{t-1}=j)=\alpha_i$.

HM-MAR model will also lead to hidden Markov model in general state space where $p$ is considered to be zero in (\ref{HM MAR}) (i.e. $Y_t$ given $Z_t$, is independent of past observations).

\section{Statistical Properties of the Model} \label{Sec StatProp}

In this section, we discuss the statistical properties of the HM-MAR model. We propose a dynamic programming approach to calculate conditional expectation and variance of the process.  We also investigate the long run behavior of the first order HM-MAR($K,1$) process, including limiting behavior of the unconditional first moment, and an appropriate upper bound for the limiting value of the variance. Finally convergence and stability of second moment is proved.

\subsection{Forecasting}

In HM-MAR model (\ref{HM MAR p1}-\ref{HM MAR}), the conditional expectation as the least square predictor (page 64 of \cite{Brockwell}) of the process $Y_t$ for $t \geq p+1$ is obtained by
\BEA
    E(Y_{t}|\mathcal{F}_{t-1})= \sum_{h=1}^K \alpha_h^{(t)}(a_{0,h}+a_{1,h}y_{t-1}+...+a_{{p},h}y_{t-{p}}), \label{cond mean}
\EEA
where $ \alpha_h^{(t)}$ is measurable $\mathcal{F}_{t-1}$.

One of the main areas for modeling conditional heteroscedasticity (changes in the conditional variance) \textbf{is the} family of ARCH models \cite{Francq}\textbf{,} originally proposed by Engle \cite{Engle} in the context of financial time series.  In the class of MTD models, MAR \cite{WL} and MAR-ARCH \cite{WL2001} models also provide a mechanism to capture this effect. However in these models only changes in conditional mean of each distribution \textbf{affect} the conditional variance of process.
The conditional variance of HM-MAR model is given by
\BEA
    Var(Y_{t}|\mathcal{F}_{t-1}) & = &
     \sum_{h=1}^K \alpha_h^{(t)}(\sigma^2_h + (a_{0,h}+ a_{1,h}y_{t-1}+...+a_{{p},h}y_{t-{p}})^2)-
    \NNL &&
    \{\sum_{h=1}^K \alpha_h^{(t)}(a_{0,h}+a_{1,h}y_{t-1}+...+a_{{p},h}y_{t-{p}})\}^2
   \NNL & = &
    \sum_{h=1}^K \alpha_h^{(t)}\sigma^2_h +  \sum_{h=1}^K \alpha_h^{(t)} \mu_{h,t}^2-
    \{\sum_{h=1}^K \alpha_h^{(t)}\mu_{h,t}\}^2
    \label{cond var}
\EEA
\textbf{in which $\mu_{h,t}=a_{0,h}+a_{1,h}y_{t-1}+...+a_{{p},h}y_{t-{p}}$ is the conditional mean of $h-$th distribution (i.e. $E[Y_t|Z_t=h,Y_1^{t-1}]$). Let $\mu_t$ be a random variable which takes values $\mu_{h,t}$ with probabilities $\alpha_h^{(t)}$ for $h=1,\cdots,K$, then $\sum_{h=1}^K \alpha_h^{(t)} \mu_{h,t}-
\{\sum_{h=1}^K \alpha_h^{(t)}\mu_{h,t}\}^2$ can be interpreted as the conditional variance of $\mu_t$ given all past observations.
 This amount is small (large) when all conditional means are equal (largely different). Relation (\ref{cond var}) shows the impact of conditional mean $\mu_{h,t}$ and weighting coefficients $ \alpha_h^{(t)}$ on the value of conditional variance of $Y_t$ given all past information}. This is the merit of the HM-MAR model and \textbf{its capability} to model conditional heteroscedasticity as a function of \textbf{simultaneous} changes in the weighting coefficients \textbf{as well as} conditional mean of each distribution.

At each time step $t$, $\alpha_h^{(t)}$ \textbf{(}in equations (\ref{cond mean}) and (\ref{cond var})\textbf{)} can be
 determined via a dynamic programming method based on forward recursion algorithm\textbf{,} proposed in remark \ref{Rem PZt|Ft-1}.

\begin{Rem} \label{Rem PZt|Ft-1}
Let $y_r^s \equiv (y_r,\cdots,y_s)$ for $s>r$, the weighting functions in the HM-MAR model (\ref{HM MAR p1}-\ref{HM MAR}) satisfy
\BEA
     \alpha_h^{(t)} & = &\frac{\sum_{m=1}^K {F(y_{p}^{t-1},Z_{t-1}=m|y_0^{p-1})\pi_{m,h}}}
                                      {\sum_{m=1}^K {F(y_{p}^{t-1},Z_{t-1}=m|y_0^{p-1})}},
                                      \label{alpha ht}
\EEA
\textbf{where} $F(y_{p}^t,z_{t}|y_1^{p-1})$ is calculated recursively as
\BEA
    F(y_{p}^t,&& \hspace{-0.29in }Z_{t}=h|y_0^{p-1})  =
    \NNL &&  \hspace{-0.29in } \sum_m F(y_{p}^{t-1},Z_{t-1}=m|y_0^{p-1})
      \pi_{m,h} \Phi(\frac{y_t-a_{0,h}-\sum_{i=1}^p a_{i,h}y_{t-i}}{\sigma_{h}}),
      \label{y_t}
\EEA
and recursion starts for $t=p$ by
\BEA
    F(y_{p},Z_{p+1}=h|y_0^{p-1}) = \rho_h
                                    \Phi(\frac{y_p-a_{0,h}-\sum_{i=1}^p a_{i,h}y_{p-i}}{\sigma_h}),
                                    \NN
\EEA
\end{Rem}

\begin{proof}
As the hidden variables $\{Z_t\}_{t \geq p}$ have Markov structure in HM-MAR model, we have \BEA
    \alpha_h^{(t)} & \hspace{-0.2in} = & \hspace{-0.2in} P(Z_t=h|y_0^{t-1}) = \sum_{m=1}^K P(Z_t=h,Z_{t-1}=m|y_0^{t-1})
                                   \NNL &=& \sum_{m=1}^K P(Z_t=h|Z_{t-1}=m,y_0^{t-1})P(Z_{t-1}=m|y_0^{t-1})
                                   \NNL &=& \sum_{m=1}^K P(Z_t=h|Z_{t-1}=m)P(Z_{t-1}=m|y_0^{t-1})
                                   \NNL &=& \frac{\sum_{m=1}^K {F(y_0^{t-1},Z_{t-1}=m)\pi_{m,h}}}
                                      {\sum_{m=1}^K {F(y_0^{t-1},Z_{t-1}=m)}}
                                   \NNL &=& \frac{\sum_{m=1}^K {F(y_{p}^{t-1},Z_{t-1}=m|y_0^{p-1})\pi_{m,h}}}
                                      {\sum_{m=1}^K {F(y_{p}^{t-1},Z_{t-1}=m|y_0^{p-1})}}, \NN
\EEA
where
\BEA
  && \hspace{-0.3in}   F(y_{p}^{t-1},Z_{t-1}=m|y_{0}^{p-1}) = \sum_{j=1}^K F(y_{p}^{t-1},Z_{t-1}=m,Z_{t-2}=j|y_{0}^{p-1}) =
                                        \NNL &&  \hspace{-0.3in}
                                      \sum_{j=1}^K\!\! F(y_{t-1}|Z_{t-1}\!=\!m,Z_{t\!-\!2}\!=\!j,y_{0}^{t\!-\!2})
                                        P(Z_{t-1}\!=\!m|Z_{t\!-\!2}\!=\!j,y_{0}^{t\!-\!2})
                                        F(y_{p}^{t\!-\!2},Z_{t\!-\!2}\!=\!j|y_{0}^{p\!-\!1})
                                        \NNL && \hspace{-0.3in}
                                        =\sum_{j=1}^K \Phi(\frac{y_{t-1}-a_{0,m}-a_{1,m}y_{t-2}-\cdots-a_{p,m}y_{t-p-1}}{\sigma_{m}})
                                        \pi_{j,m}
                                        F(y_{p}^{t-2},Z_{t-2}=j|y_{0}^{p-1}),
                                        \NN
\EEA
in which the last equality implies by (\ref{HM MAR}) and the recursion begins with (\ref{HM MAR p1}).
\end{proof}

Another characteristic of HM-MAR is modeling the \textbf{all} past observations and  \textbf{benefits} from a dynamic programming approach\textbf{. This will in turn} minimize the volume of calculations  for forecasting. \textbf{T}he intermediate results \textbf{and in} fact the last state  $ F(y_{p+1},...,y_t,Z_{t}=h|y_1,...,y_p)$ \textbf{is stored} for different values of $Z_t$ \textbf{which could be used to update the process}, see (\ref{alpha ht}-\ref{y_t}).

\subsection{Stability}
In this \textbf{section}, we investigate the stability of moments for the nonlinear process $\{Y_t\}_{t=0}^\infty$ that admits a HM-MAR($K,1$) model.
This process \textbf{is} represented as a random coefficient autoregressive process of order one,
in which the autoregressive coefficients are functions of the latent random variables ,$\{Z_t\}_{t \geq 1}$, \textbf{(see Equations} (\ref{Initial Probabilities})-(\ref{Mu})\textbf{)}.
Let random variable\textbf{s} \textbf{and $\sigma_{Z_{t}}$ respectively} $a_{i,Z_t}$ \textbf{take} values $\{a_{i,1}\,,\cdots,\,a_{i,K}\}$ for $i=0,1$, and $\{\sigma_{1}\,,\cdots,\,\sigma_{K}\}$, \textbf{where} $a_{i,j}$ and $\sigma_j$, $j=1,\cdots,K$ are used in HM-MAR model  (\ref{HM MAR p1}-\ref{HM MAR}) with $p=1$.
 \textbf{We} consider
\BEA
    Y_t= a_{0,Z_{t}} + a_{1,Z_{t}}Y_{t-1} + \sigma_{Z_{t}} \varepsilon_{t}, \label{HM MAR1}
\EEA
where $\{\varepsilon_t\}_{t \geq 1}$ is a Gaussian IID(0,1) process, independent of the hidden process $\{Z_t\}_{t \geq 1}$. The conditional distribution of the process $Y_t$ in \textbf{E}quation (\ref{HM MAR1}) is determined as
\BEA
    F(y_t|\mathcal{F}_{t-1}) = \sum_{h=1}^K  P(Z_{t}=h|\mathcal{F}_{t-1}) F(y_t|Z_{t}=h,\mathcal{F}_{t-1}), \NN
\EEA
in which $P(z_{t}=h|\mathcal{F}_{t-1})= \alpha_t^h$ is given by remark \ref{Rem PZt|Ft-1}\textbf{. By the} Gaussian distribution of $\varepsilon_t$ in (\ref{HM MAR1}), we have
\BEA
    F(y_t|Z_{t}=h,\mathcal{F}_{t-1}) & = &
                                \Phi(\frac{y_t-a_{0,h}-a_{1,h}y_{t-1}}{\sigma_{h}}). \NN
\EEA
Thus (\ref{HM MAR p1}-\ref{HM MAR}) implies that $\{Y_t\}_{t=0}^\infty$  admits HM-MAR($K,1$) model.

\textbf{Notice} that the process $\{Y_t\}_{t=0}^\infty$ is not necessarily a  Markov process, however the extended process $X=\{X_t\}_{t=1}^\infty$ with $X_t=(Z_t,\bar{Y}_t=(Y_t,Y_{t-1},...,Y_{t-p})^\prime)^\prime$ is Markov \cite{Yao}.

Timmermann \cite{Timmermann} derived the moments of a class of stationary Markov switching models with state-dependent
autoregressive dynamics and conditional mean, $\mu_{Z_t}$.
Our approach for \textbf{deriving} the limiting behavior of first and second moments of the process $Y_t$ is not based on the \textbf{stationary} assumption of \textbf{the } model.

\noindent
Let's define the $K \times K$ diagonal matrixes
    \BEA
        \boldsymbol{\phi}_i & = & diag(a_{i,1},\, \cdots ,\, a_{i,K}), \qquad i=0,1, \NNL
        \boldsymbol{\sigma} & = & diag(\sigma_{1},\, \cdots ,\, \sigma_{K}), \NN
    \EEA
for possible values of random variables $a_{i,Z_t}$ and $\sigma_{Z_{t}}$ in equation (\ref{HM MAR1}) \textbf{where} $\mathbf{1}=(1,\cdots,1)^\prime$ \textbf{is} a $K \times 1$ vector.

\begin{Lem} \label{Lem EphiZn|Z1}
    Let $\{Y_t\}_{t=0}^\infty$ be a HM-MAR($K,1$) process defined by (\ref{HM MAR1}), then for $n \geq 2$
    \BEA
        \left(\begin{array}{c}E[\prod_{t=2}^n a_{1,Z_{t}}|Z_{1}=1] \\
                                \vdots \\
                              E[\prod_{t=2}^n a_{1,Z_{t}}|Z_{1}=K]\end{array}\right)
                              = (P \boldsymbol{\phi}_1)^{n-1} \mathbf{1}.
                     \NN
    \EEA
\end{Lem}

\begin{proof}
    By the Markov property of  $\{Z_t\}_{t=1}^\infty$ we have that $$E[a_{1,Z_t}|\sigma\{Z_{s},s \leq t-1\}]=E[a_{1,Z_t}|Z_{t-1}].$$ So
    \BEA
        \NewAND E[\prod_{t=2}^{n} a_{1,Z_{t}}|Z_{1}=k]  =
            \sum_{Z_{2},\cdots,Z_n} (\prod_{t=2}^{n} a_{1,Z_{t}}) P(Z_{2},\cdots,Z_n|Z_{1}=k)
        \NNL & = & \sum_{Z_{2},\cdots,Z_n} (\prod_{t=2}^{n} a_{1,Z_{t}})
        P(Z_{3},\cdots,Z_n|Z_{1},Z_{1}=k) P(Z_{2}|Z_{1}=k)
        \NNL & = &
        \sum_{Z_{2}}\{\sum_{Z_{3},\cdots,Z_n} (\prod_{t=3}^{n} a_{1,Z_{t}})  P(Z_{3},\cdots,Z_n|Z_{2}) \}
        a_{1,Z_{2}} P(Z_{2}|Z_{1}=k)
        \NNL & = & E[E[\prod_{t=3}^{n} a_{1,Z_{t}}|Z_{2}]a_{1,Z_{2}}|Z_{1}=k]. \NN
    \EEA
    So for vector of conditional expectations of $\prod_{t=2}^{n} a_{1,Z_{t}}$ given different values of $Z_{1}$, we have the following recursive equation
    \BEA
       \NewAND \hspace{-0.1in} \left(\begin{array}{c}E[\prod_{t=2}^{n} a_{1,Z_{t}}|Z_{1}=1] \\
                                \vdots \\
                              E[\prod_{t=2}^{k+1} a_{1,Z_{t}}|Z_{1}=K]\end{array}\right)
       = \left(\begin{array}{c}E[E[\prod_{t=2}^{n} a_{1,Z_{t}}|Z_{2}]|Z_{1}=1] \\
                                \vdots \\
                              E[E[\prod_{t=2}^{n} a_{1,Z_{t}}|Z_{2}]|Z_{1}=K]\end{array}\right)
       \NNL \NewAND                      =
       \left(\begin{array}{c}E[E[\prod_{t=3}^{n} a_{1,Z_{t}}|Z_{2}]a_{1,Z_{2}}|Z_{1}=1] \\
                                \vdots \\
                              E[E[\prod_{t=3}^{n} a_{1,Z_{t}}|Z_{2}]a_{1,Z_{2}}|Z_{1}=K]\end{array}\right)
       \NNL \NewAND =
       \left(\begin{array}{c}\sum_{i=1}^K E[\prod_{t=3}^{n} a_{1,Z_{t}}|Z_{2}=i]a_{1,i} \pi_{1i} \\
                                \vdots \\
                              \sum_{i=1}^K E[\prod_{t=3}^{n} a_{1,Z_{t}}|Z_{2}=i]a_{1,i} \pi_{Ki}
                              \end{array}\right)
       \NNL\NewAND =
       \left(\begin{array}{ccc} \pi_{11} & \cdots & \pi_{1K} \\
                                \vdots & \vdots & \vdots \\
                                \pi_{K1} & \cdots & \pi_{KK} \end{array}\right)
              \left(\begin{array}{ccc} a_{1,1} & 0 & \cdots \\
                                \vdots & \vdots & \vdots \\
                                0 & \cdots & a_{1,K} \end{array}\right)
              \left(\begin{array}{c}E[\prod_{t=3}^{n} a_{1,Z_{t}}|Z_{2}=1] \\
                                \vdots \\
                              E[\prod_{t=3}^{n} a_{1,Z_{t}}|Z_{2}=K]\end{array}\right)
        \NNL \NewAND =
                P \boldsymbol{\phi}_1
                \left(\begin{array}{c}E[\prod_{t=3}^{n} a_{1,Z_{t}}|Z_{2}=1] \\
                                \vdots \\
                              E[\prod_{t=3}^{n} a_{1,Z_{t}}|Z_{2}=K]\end{array}\right), \label{Recursive E}
\EEA
in which the recursion starts at $t=n-1$ as
    \BEA
        \NewAND \left(  \begin{array}{c} E[a_{1,Z_{n}}|Z_{n-1}=1] \\
                                \vdots \\
                              E[a_{1,Z_{n}}|Z_{n-1}=K]\end{array}  \right)
        = \left(\begin{array}{c}\sum_{i=1}^K a_{1,i} \pi_{1k} \\
                                                            \vdots \\
                                                            \sum_{i=1}^K a_{1,i} \pi_{Kk}\end{array}\right)
                                \NNL &=& \left(\begin{array}{ccc} \pi_{11}a_{1,1} & \cdots & \pi_{1K}a_{1,K} \\
                                                            \vdots & \vdots & \vdots \\
                                                            \pi_{K1}a_{1,1} & \cdots & \pi_{KK}a_{1,K} \end{array}\right)
                                                            \left(\begin{array}{c} 1 \\
                                                            \vdots \\
                                                            1
                                                            \end{array}\right)
                                \NNL &=& \left(\begin{array}{ccc} \pi_{11} & \cdots & \pi_{1K} \\
                                                            \vdots & \vdots & \vdots \\
                                                            \pi_{K1} & \cdots & \pi_{KK} \end{array}\right)
                                          \left(\begin{array}{ccc} a_{1,1} & 0 & \cdots \\
                                                            \vdots & \vdots & \vdots \\
                                                            0 & \cdots & a_{1,K} \end{array}\right)
                                                            \left(\begin{array}{c} 1 \\
                                                            \vdots \\
                                                            1
                                                            \end{array}\right)
                \nonumber \\ & = & P \boldsymbol{\phi}_1 \mathbf{1}. \NN
    \EEA
Thus the solution of recursive equation (\ref{Recursive E}) is given by
    \BEA
       \left(\begin{array}{c}E[\prod_{t=2}^n a_{1,Z_{t}}|Z_{1}=1] \\
                                \vdots \\
                              E[\prod_{t=2}^n a_{1,Z_{t}}|Z_{1}=K]\end{array}\right)  = (P \boldsymbol{\phi}_1)^{n-1} \mathbf{1}. \NN
    \EEA
\end{proof}
\begin{Lem} \label{Lem EphiZn}
Let $\{Z_t\}_{t=1}^\infty$ be a Markov chain starting with invariant probability measure $\boldsymbol{\mu}$ defined by (\ref{Mu}), then under conditions of the lemma  \ref{Lem EphiZn|Z1}
 \BEA
    E[\prod_{t=2}^n a_{1,Z_t}a_{0,Z_{1}}]=\boldsymbol{\mu}^\prime \boldsymbol{\phi}_0
                                                    (P \boldsymbol{\phi}_1)^{n-1} \mathbf{1}. \NN
 \EEA
\end{Lem}
\begin{proof}
By lemma \ref{Lem EphiZn|Z1}, we have
\BEA
    && E[\prod_{t=2}^n a_{1,Z_t}a_{0,Z_{1}}]  =  E[E[\prod_{t=2}^n a_{1,Z_t}|Z_{1}]a_{0,Z_{1}}]
    \NNL &  & =
    \sum_{k=1}^K \alpha_k a_{k,0}E[\prod_{t=2}^n a_{1,Z_t}|Z_{p}=k]
    \NNL &  & =(\alpha_1\,,\cdots,\,\alpha_K)^\prime
        \left[\begin{array}{cccc} a_{0,1} & 0 & \cdots & 0 \\
                                                        0 & a_{0,2} & \cdots & 0 \\
                                                        \vdots &  & \cdots & \vdots \\
                                                         0 & 0 & \cdots & a_{0,K} \end{array}\right]
    \left[\begin{array}{c}E[\prod_{t=2}^n a_{1,Z_t}|Z_{p}=1] \\
                                \vdots \\
                          E[\prod_{t=2}^n a_{1,Z_t}|Z_{p}=K]\end{array}\right]
    \NNL &  & = \boldsymbol{\mu}^\prime \boldsymbol{\phi}_0(P \boldsymbol{\phi}_1)^{n-1} \mathbf{1}.
\EEA
\end{proof}

\begin{Lem} \label{Lem Datta}
If all eigenvalues of $P \boldsymbol{\phi}_1$ lie inside the unite circle then under conditions of lemma \ref{Lem EphiZn}
    \begin{enumerate}[label=\roman{*}., ref=(\roman{*})]
        \item $ \lim_{m \rightarrow \infty} E[\prod_{n=2}^{m+1} a_{1,Z_n}a_{0,Z_{1}}] = 0,$

        \item $lim_{t \rightarrow \infty} \sum_{m=0}^t E[\prod_{n=2}^{m+1} a_{1,Z_n}a_{0,Z_{1}}] = \boldsymbol{\mu}^\prime \boldsymbol{\phi}_0 (I- P\boldsymbol{\phi}_1)^{-1}\mathrm{1}. $
    \end{enumerate}
Also if all eigenvalues of $P \boldsymbol{\phi}_1^2$ lie inside the unite circle then
    \begin{enumerate}[label=\roman{*}., ref=(\roman{*})]
        \item $ \lim_{m \rightarrow \infty} E[(\prod_{n=2}^{m+1} a_{1,Z_n}a_{0,Z_{1}})^2] = 0,$
        \item $ \lim_{t \rightarrow \infty} \sum_{m=0}^t E[(\prod_{n=2}^{m+1} a_{1,Z_n}a_{0,Z_{1}})^2] = \boldsymbol{\mu}^\prime \boldsymbol{\phi}_0^2 (I- P\boldsymbol{\phi}_1^2)^{-1}\mathrm{1}.$
    \end{enumerate}
\end{Lem}
\begin{proof}
The first part is an immediate result of lemma \ref{Lem EphiZn} and Datta (page 508 of \cite{Datta2004}) and for the second part:
\BEA
    \lim_{t \rightarrow \infty} \sum_{m=0}^t E[\prod_{n=2}^{m+1} a_{1,Z_n}a_{0,Z_{1}}] & = &
     \lim_{t \rightarrow \infty} \sum_{m=0}^t \boldsymbol{\mu}^\prime \boldsymbol{\phi}_0
                                                    (P \boldsymbol{\phi}_1)^{m} \mathbf{1}
    \NNL & = &  \boldsymbol{\mu}^\prime \boldsymbol{\phi}_0 (I- P\boldsymbol{\phi}_1)^{-1}\mathrm{1}, \label{lim E sum prod a_i}
\EEA
in which the last equality holds by Datta (page 511 of \cite{Datta2004}). The rest of proof can be done in a similar way by conducting a result similar to lemma \ref{Lem EphiZn} as $E[(\prod_{n=2}^{m+1} a_{1,Z_n}a_{0,Z_{1}})^2]=\boldsymbol{\mu}^\prime \boldsymbol{\phi}_0^2 (P\boldsymbol{\phi}_1^2)^{m}\mathrm{1}$.
\end{proof}

\begin{Lem} \label{Lem lim E prod phi Y0}
If $E[Y_0^2] < \infty$ then under conditions of lemma \ref{Lem Datta}
   $$\lim_{t \rightarrow \infty} E[\prod_{i=1}^{t} a_{1,Z_{i}}Y_{0}] = 0.$$
\end{Lem}
\begin{proof}
By Cauchy Schwarz inequality we have
\BEA
    [Cov(\prod_{i=1}^{t} a_{1,Z_{i}},Y_{0})]^2 < Var(\prod_{i=1}^{t} a_{1,Z_{i}})Var(Y_{0}), \NN
\EEA
by lemma \ref{Lem Datta} we can deduce that
\BEA
    \lim_{t \rightarrow \infty} Var(\prod_{i=1}^{t} a_{1,Z_{i}}) & = & 0, \NN
\EEA
and since $Var(Y_0) < \infty$, so
\BEA
    \lim_{t \rightarrow \infty} Cov(\prod_{i=1}^{t} a_{1,Z_{i}},Y_{0}) = 0, \NN
\EEA
thus
\BEA
    \lim_{t \rightarrow \infty} E[\prod_{i=1}^{t} a_{1,Z_{i}}Y_{0}] = \lim_{t \rightarrow \infty} E[\prod_{i=1}^{t} a_{1,Z_{i}}]E[Y_{0}]=0, \NN
\EEA
in which the last equality can be verified by lemma \ref{Lem Datta} and the fact that $E[Y_0]$ is finite by the assumption that $E[Y^2_0] < \infty$ (page 274 of \cite{Billingsley}).
\end{proof}

\begin{Theo} \label{Theo Mean HM-MAR1}
    Let $\{Y_t\}_{t=0}^\infty$ follows the HM-MAR($K,1$) model, defined by (\ref{HM MAR1}), and the  following assumptions hold
    \begin{enumerate}[label=\roman{*}., ref=(\roman{*})]
        \item $\{Z_t\}_{t>1}$ is an ergodic Markov chain  starting from its invariant probability measure $\boldsymbol{\mu}$ specified in equation (\ref{Mu}), \label{Cond ergodic Theo Mean HM-MAR1}
        \item $E[Y_0^2] < \infty$, \label{Cond finite mean Theo Mean HM-MAR1}
        \item All eigenvalues of $P \boldsymbol{\phi}_1$ and $P \boldsymbol{\phi}_1^2 $ lie inside the unit circle,
        \label{Cond eigenvalue Theo Mean HM-MAR1}
    \end{enumerate}
    then the process is asymptotically stable in mean and
    \BEA
        \lim_{t \rightarrow \infty} E[Y_t]= \boldsymbol{\mu} \boldsymbol{\phi}_0 (I -  P \boldsymbol{\phi}_1)^{-1}\mathbf{1}.  \label{lim E Yt}
    \EEA
\end{Theo}
\begin{proof}
    Iterating equation (\ref{HM MAR1}), we get
    \BEA
      \hspace{-0.5in}
       Y_t  &=&   a_{0,Z_{t}} + a_{1,Z_{t}}Y_{t-1} + \sigma_{Z_{t}} \varepsilon_{t}  \nonumber \\
            &=&    a_{0,Z_{t}} + a_{1,Z_{t}} a_{0,Z_{t-1}} + a_{1,Z_{t}} \sigma_{Z_{t-1}} \varepsilon_{t-1}  + \sigma_{Z_{t}} \varepsilon_{t} + a_{1,Z_{t}}a_{1,Z_{t-1}}Y_{t-2}  \nonumber \\
            &=& \hspace{-0.1in} \sum_{m=0}^{t-1} \prod_{i=0}^{m-1} a_{1,Z_{t-i}}(a_{0,Z_{t-m}} +
             \sigma_{Z_{t-m}} \varepsilon_{t-m}) +
            \prod_{i=0}^{t-1} a_{1,Z_{t-i}}Y_{0}. \label{expanded Y for t-i}
    \EEA
    Let $u=t-i$ in (\ref{expanded Y for t-i}) to get
    \BEA
        Y_t  &=&   \sum_{m=0}^{t-1} \prod_{u=t-m+1}^{t} a_{1,Z_{u}}(a_{0,Z_{t-m}} +
             \sigma_{Z_{t-m}} \varepsilon_{t-m}) +
            \prod_{u=1}^{t} a_{1,Z_{u}}Y_{0}
            \NNL &=&
            \sum_{m=0}^{t-1} \prod_{u=2}^{m+1} a_{1,Z_{u}}(a_{0,Z_{1}} +
             \sigma_{Z_{1}} \varepsilon_{t-m}) +
            \prod_{u=1}^{t} a_{1,Z_{u}}Y_{0},  \label{expanded Y}
    \EEA
    where the last equality follows from the strict stationarity property of $\{Z_t\}_{t=1}^\infty$ (page 35 of \cite{Fan}), which implies by assumption \ref{Cond ergodic Theo Mean HM-MAR1} of theorem.
    Also by the independence assumption of $\{\varepsilon_t\}$ from $\{Z_t\}_{t=1}^\infty$ in (\ref{HM MAR1}):
    \BEA
      \hspace{-0.4in} \lim_{t \rightarrow \infty} E[\sum_{m=0}^{t-1} \prod_{u=2}^{t-m} a_{1,Z_{u}} \sigma_{Z_{1}} \varepsilon_{t-m}]
       & = &
         \lim_{t \rightarrow \infty} E[\sum_{m=0}^{t-1} \prod_{u=2}^{m+1} a_{1,Z_{u}}]E[\varepsilon_{t-m}] = 0.
         \label{E sum prod phi sigma}
    \EEA
    Thus by lemma  \ref{Lem lim E prod phi Y0} and (\ref{expanded Y}- \ref{E sum prod phi sigma}) we have that
    \BEA
       \lim_{t \rightarrow \infty }E[Y_t] = \lim_{t \rightarrow \infty } E[\sum_{m=0}^{t-1} \prod_{u=2}^{m+1} a_{1,Z_{u}}a_{0,Z_{1}}], \NN
    \EEA
    so by assumption \ref{Cond eigenvalue Theo Mean HM-MAR1} and lemma \ref{Lem Datta}, we get (\ref{lim E Yt}).
\end{proof}

One interesting feature of \textbf{T}heorem \ref{Theo Mean HM-MAR1} is that HM-MAR model \textbf{could consist of some}  explosive (with $a_1 \geq 1$) and non-explosive autoregressive processes and it remains asymptotically stable in mean.

\begin{Def} \label{Def spectral raduis}
Let $\lambda$ be the spectral radius of
$$A \equiv \mathbf{1}(P\boldsymbol{\phi}_1^2\mathbf{1})^\prime \boldsymbol{I} = diag(E[a_{1,Z_{t}}^2|Z_{t-1}=1],\cdots,E[a_{1,Z_{t}}^2|Z_{t-1}=K]).$$
\end{Def}

\begin{Lem} \label{Lem E (sum phi1 phi0)2}
Let spectral radius $\lambda$ to be as in definition \ref{Def spectral raduis}. If $\lambda$ lies inside the unit circle then under conditions of lemma \ref{Lem EphiZn}
\BEA
   && \lim_{t \rightarrow \infty}  E[(\sum_{m=0}^{t-1} \prod_{i=2}^{m+1}  a_{1,Z_{i}} a_{0,Z_{1}})^2] <
    2(\frac{1+\boldsymbol{\mu}^\prime \boldsymbol{\phi}_0^2 \mathbf{1}}{1-\lambda^{1/2}})^2
    <  \infty. \label{ineq E (sum prod a_i)^2}
\EEA
\textbf{Furthermore if} $E[Y_0^{2+\epsilon}] < \infty ,\,\, \epsilon > 0$ then
\BEA
    \lim_{t \rightarrow \infty} E[\prod_{i=1}^{t} a_{1,Z_{i}}^2 Y_0^2 ] =0. \NN
\EEA
\end{Lem}
\begin{proof}
By definition of spectral radius wee have that the absolute values of all eigenvalues of $A$ are less than or equal to $\lambda$, so
by the lemma assumption about $\lambda$, we have that $E[a_{1,Z_t}^2|Z_{t-1}=k] \leq \lambda <1$ for all values of $k=1,\cdots,K$, thus by the method of iterative conditioning
\BEA
    E[\prod_{u=2}^{m+1} a_{1,Z_{u}}^2a_{0,Z_{1}}^2] &=& E[E[\prod_{u=2}^{m+1} a_{1,Z_{u}}^2a_{0,Z_{1}}^2|\sigma\{Z_1^m\}]]
    \NNL & =& E[E[a_{1,Z_{m+1}}^2|\sigma\{Z_1^m\}]\prod_{u=2}^{m} a_{1,Z_{u}}^2a_{0,Z_{1}}^2]
    \NNL & \leq & \lambda E[\prod_{u=2}^{m} a_{1,Z_{u}}^2a_{0,Z_{1}}^2],
    \label{max E prod phi1 phi0}
\EEA
in which $\sigma\{Z_1^m\} \equiv \sigma\{Z_1,\cdots,Z_m\}$.  Iterating (\ref{max E prod phi1 phi0}) we get
\BEA
    E[\prod_{u=2}^{m+1} a_{1,Z_{u}}^2a_{0,Z_{1}}^2]
    \leq \lambda^{m} E[a_{0,Z_{1}}^2]
    = \lambda^{m} \boldsymbol{\mu}^\prime \boldsymbol{\phi}_0^2 \mathbf{1},  \label{ineq E prod phi1 phi0}
\EEA
thus
\BEA
     \lim_{t \rightarrow \infty} \sum_{m=0}^{t-1} E[\prod_{u=2}^{m+1} a_{1,Z_{u}}^2a_{0,Z_{1}}^2]
     \leq
     \boldsymbol{\mu}^\prime \boldsymbol{\phi}_0^2 \mathbf{1}
     (\lim_{t \rightarrow \infty} \sum_{m=0}^{t-1} \lambda^{m})
     = \frac{\boldsymbol{\mu}^\prime \boldsymbol{\phi}_0^2 \mathbf{1}}{1-\lambda} . \label{ineq E prod phi^2}
\EEA
Now by Cauchy Schwarz inequality,
\BEA
    \hspace{-0.3in} E^2[(\prod_{i=2}^{m+1}a_{1,Z_{i}}a_{0,Z_{1}})
    && \hspace{-0.25in}(\prod_{j=0}^{n+1}a_{1,Z_{j}}a_{0,Z_{1}})]
    \NNL && \leq
    E[\prod_{i=2}^{m+1}a_{1,Z_{i}}^2a_{0,Z_{1}}^2]
    E[\prod_{j=2}^{n+1}a_{1,Z_{j}}^2a_{0,Z_{1}}^2], \NN
\EEA
thus
\BEA
    E[(\prod_{i=2}^{m+1}a_{1,Z_{i}}a_{0,Z_{1}})(\prod_{j=2}^{n+1}a_{1,Z_{j}}a_{0,Z_{1}})]
    \leq \boldsymbol{\mu}^\prime \boldsymbol{\phi}_0^2 \mathbf{1}\lambda^{(m+n)/2}, \NN
\EEA
and summing up for different values of $m \neq n = 0$ to $\infty$,
\BEA
    &&\hspace{-0.3in} \sum_{m \neq n=0}^{\infty} E[(\prod_{i=2}^{m+1}a_{1,Z_{i}}a_{0,Z_{1}})
    (\prod_{j=2}^{n+1}a_{1,Z_{j}}a_{0,Z_{1}})]
    <
    \sum_{m \neq n=0}^{\infty} \boldsymbol{\mu}^\prime \boldsymbol{\phi}_0 \mathbf{1}\lambda^{(m+n)/2}
    \NNL && <
    (\sum_{m =0}^{\infty} \boldsymbol{\mu}^\prime \boldsymbol{\phi}_0^2 \mathbf{1}\lambda^{(m)/2})^2
    = (\frac{\boldsymbol{\mu}^\prime \boldsymbol{\phi}_0^2 \mathbf{1}}{1-\lambda^{1/2}})^2. \label{ineq sum E2 prod phi}
\EEA
Now by (\ref{ineq E prod phi^2}) and (\ref{ineq sum E2 prod phi}) we have
\BEA
    && \hspace{-0.3in}
    \lim_{t \rightarrow \infty}
     E[(\sum_{m=0}^{t-1} \prod_{i=2}^{m+1}a_{1,Z_{i}} a_{0,Z_{1}})^2]=
    \NNL && \hspace{-0.3in}
     \sum_{m=0}^{\infty} E[(\prod_{i=2}^{m+1}a_{1,Z_{i}} a_{0,Z_{1}})^2] +
      2 \sum_{m \neq n=0}^{\infty} E[(\prod_{i=2}^{m+1}a_{1,Z_{i}}a_{0,Z_{1}})
    (\prod_{j=2}^{n+1}a_{1,Z_{j}}a_{0,Z_{1}})]
    \NNL && \hspace{-0.3in} <
    \frac{\boldsymbol{\mu}^\prime \boldsymbol{\phi}_0^2 \mathbf{1}}{1-\lambda}
    +
    2(\frac{\boldsymbol{\mu}^\prime \boldsymbol{\phi}_0^2 \mathbf{1}}{1-\lambda^{1/2}})^2
    <
    2(\frac{1+\boldsymbol{\mu}^\prime \boldsymbol{\phi}_0^2 \mathbf{1}}{1-\lambda^{1/2}})^2. \NN
\EEA

\noindent
Now by Holder inequality (page 80 of \cite{Billingsley}),
$$E[ a_{1,Z_{1}}^2 Y_0^2] < E^{1/u}[a_{1,Z_{1}}^{2u}]E^{1/v}[Y_0^{2v}] = (\boldsymbol{\mu}^\prime \boldsymbol{\phi}_1^{2u}\mathbf{1})^{1/u}E^{1/v}[Y_0^{2v}] < \infty,$$
in which $u,v >1$ and $1/u + 1/v =1$, so for $v=1+\epsilon/2$ we set $u=v/(v-1)$ and $(\boldsymbol{\mu}^\prime \boldsymbol{\phi}_1^{2u}\mathbf{1})^{1/u} < \infty$.
\noindent
Thus by inequality (\ref{ineq E prod phi1 phi0}) and the fact that $\lambda < 1$, we have
\BEA
   \lim_{t \rightarrow \infty} E[\prod_{i=1}^{t} a_{1,Z_{i}}^2 Y_0^2] =  \lim_{t \rightarrow \infty} \lambda^{t-1} E[ a_{1,Z_{1}}^2 Y_0^2] = 0. \NN
\EEA
\end{proof}

\noindent
Thus by lemmas \ref{Lem Datta} and \ref{Lem E (sum phi1 phi0)2} \textbf{, we got the following inequality}
\BEA
    \lim_{t \rightarrow \infty}
    Var(\sum_{m=0}^{t-1} \prod_{i=2}^{m+1}
    a_{1,Z_{i}} a_{0,Z_{1}})
     \hspace{-0.05in}< \hspace{-0.05in} 2(\frac{1+\boldsymbol{\mu}^\prime \boldsymbol{\phi}_0^2 \mathbf{1}}{1-\lambda^{1/2}})^2
    \hspace{-0.05in} - \hspace{-0.05in}
    (\boldsymbol{\mu} \boldsymbol{\phi}_0 (I- P\boldsymbol{\phi}_1)^{-1}\mathbf{1})^2 \label{LimVar1}.
\EEA

\begin{Theo} \label{Theo second HM-MAR1}
    Let $\{Y_t\}_{t=0}^\infty$ \textbf{follow} the HM-MAR($K,1$) model defined by (\ref{HM MAR1}) \textbf{with} $\lambda$ \textbf{as in }definition \ref{Def spectral raduis}.
    If \textbf{the} conditions of theorem \ref{Theo Mean HM-MAR1} \textbf{hold} and
    \begin{enumerate}[label=\roman{*}., ref=(\roman{*})]
        \item $E[Y_0^{2+\epsilon}] < \infty,\qquad \epsilon > 0,$
        \item $\lambda<1,$
    \end{enumerate}
    then the process has finite second moment and
    \BEA
        \lim_{t \rightarrow \infty} E(Y_t^2) \leq
          2(\frac{1+\boldsymbol{\mu}^\prime \boldsymbol{\phi}_0^2 \mathbf{1}}{1-\lambda^{1/2}})^2
          +  \boldsymbol{\mu} \boldsymbol{\sigma}^2 (I- P\boldsymbol{\phi}_1^2)^{-1}\mathbf{1}. \label{lim E Y2 ineq}
    \EEA
\end{Theo}
\begin{proof}
Using (\ref{expanded Y}) we have
\BEA
            E[Y^2_t] &=& E[\{\sum_{m=0}^{t-1} \prod_{i=2}^{m+1} a_{1,Z_{i}} a_{0,Z_{1}}\}^2] +
                    \NN  E[\{\sum_{m=0}^{t-1} \prod_{i=2}^{m+1} a_{1,Z_{i}} \sigma_{Z_{1}}\varepsilon_{t-m}\}^2] +
                    \NNL && \hspace{-0.4in}
                    E[\prod_{i=1}^{t} a_{1,Z_{i}}^2Y_{0}^2] +
                    2E[\{\sum_{m=0}^{t-1} \prod_{i=2}^{m+1} a_{1,Z_{i}} a_{0,Z_{1}}\} \prod_{i=1}^{t} a_{1,Z_{i}}Y_{0}] +
                    \NNL &&  \hspace{-0.4in}
                    2E[(\sum_{m=0}^{t-1} \prod_{i=2}^{m+1} a_{1,Z_{i}} a_{0,Z_{1}}+\prod_{i=1}^{t} a_{1,Z_{i}}Y_{0})
                    (\sum_{m=0}^{t-1} \prod_{i=2}^{m+1} a_{1,Z_{i}} \sigma_{Z_{1}}\varepsilon_{t-m})], \label{E Y2}
\EEA
by independence of Gaussian IID(0,1) process, $\{\varepsilon_t\}$ from $\{Z_t\}$,  (as indicated in (\ref{HM MAR1})), we have
\BEA
    E[(\sum_{m=0}^{t-1} \prod_{i=2}^{m+1} a_{1,Z_{i}} a_{0,Z_{1}})(\sum_{m=0}^{t-1} \prod_{i=2}^{m+1} a_{1,Z_{i}} \sigma_{Z_{1}}\varepsilon_{t-m})]=0.
     \label{E sum prod phi varepsilon}
\EEA

\noindent
Also by Cauchy Schwarz inequality we have that
\BEA
    [Cov(\{\sum_{m=0}^{t-1} \prod_{i=2}^{m+1} \NewAND a_{1,Z_{i}} a_{0,Z_{1}}\} ,\prod_{i=1}^{t} a_{1,Z_{i}}Y_0)]^2
    \NNL && \leq
    Var(\sum_{m=0}^{t-1} \prod_{i=2}^{m+1} a_{1,Z_{i}} a_{0,Z_{1}})
    Var(\prod_{i=1}^{t} a_{1,Z_{i}}Y_0) \NN,
\EEA
lemmas \ref{Lem lim E prod phi Y0} and \ref{Lem E (sum phi1 phi0)2} imply that $\lim_{t \rightarrow \infty} Var(\prod_{i=1}^{t} a_{1,Z_{i}}Y_0)=0$, so by (\ref{LimVar1}) we have
\BEA
    && \lim_{t \rightarrow \infty} [Cov(\{\sum_{m=0}^{t-1} \prod_{i=2}^{m+1} a_{1,Z_{i}} a_{0,Z_{1}}\} ,\prod_{i=1}^{t} a_{1,Z_{i}}Y_0)]^2 = 0, \NN
\EEA
so we get
\BEA
    && \hspace{-0.5in}
    \lim_{t \rightarrow \infty} E[\{\sum_{m=0}^{t-1} \prod_{i=2}^{m+1} a_{1,Z_{i}} a_{0,Z_{1}}\}\prod_{i=1}^{t} a_{1,Z_{i}}Y_0]
    \NNL && =
    \lim_{t \rightarrow \infty}
    E[\{\sum_{m=0}^{t-1} \prod_{i=2}^{m+1} a_{1,Z_{i}} a_{0,Z_{1}}\}]E[\prod_{i=1}^{t} a_{1,Z_{i}}Y_0]
    =0,
    \label{lim E prod phi1 prod phi0 Y0}
\EEA
in which the last equality follows by lemma \ref{Lem Datta} and lemma \ref{Lem lim E prod phi Y0}.
By a similar method as for (\ref{lim E prod phi1 prod phi0 Y0}) we get
\BEA
   \lim_{t \rightarrow \infty} E[&& \hspace{-0.25in}(\{\sum_{m=0}^{t-1} \prod_{i=2}^{m+1}a_{1,Z_{i}} \sigma_{0,Z_{1}}\varepsilon_{t-m}\})(\prod_{i=1}^{t} a_{1,Z_{i}}Y_0)]
   \NNL &&  = \lim_{t \rightarrow \infty} E[\{\sum_{m=0}^{t-1} \prod_{i=2}^{m+1} a_{1,Z_{i}} a_{0,Z_{1}}\varepsilon_{t-m}\}]E[\prod_{i=1}^{t} a_{1,Z_{i}}Y_0] = 0.
   \label{E sum prod phi varepsilon Y0}
\EEA
Thus collecting results, by lemma \ref{Lem lim E prod phi Y0}, (\ref{E Y2}-\ref{E sum prod phi varepsilon Y0}) we have
\BEA
   \lim_{t \rightarrow \infty}  E[Y^2_t] =
    \lim_{t \rightarrow \infty} \{ E[\{\sum_{m=0}^{t-1} \prod_{i=2}^{m+1} a_{1,Z_{i}} a_{0,Z_{1}}\}^2] +
          E[\sum_{m=0}^{t-1} \{\prod_{i=2}^{m+1} a_{1,Z_{i}} \sigma_{Z_{1}}\}^2]\}. \NN
\EEA
Now by lemma \ref{Lem Datta},
$\lim_{t \rightarrow \infty} E[\sum_{m=0}^{t-1} \{\prod_{i=2}^{m+1} a_{1,Z_{i}} \sigma_{Z_{1}}\}^2] = \boldsymbol{\mu} \boldsymbol{\sigma}^2 (I- P\boldsymbol{\phi}_1^2)^{-1}\mathbf{1}$, so using lemma \ref{Lem E (sum phi1 phi0)2} we get (\ref{lim E Y2 ineq}).
\end{proof}

\begin{Rem}
An immediate consequence of \textbf{theorems} \ref{Theo Mean HM-MAR1} and \ref{Theo second HM-MAR1} is that
\BEA
        \lim_{t \rightarrow \infty} Var(Y_t) \leq
           2(\frac{1+\boldsymbol{\mu}^\prime \boldsymbol{\phi}_0^2 \mathbf{1}}{1-\lambda^{1/2}})^2
          +  \boldsymbol{\mu} \boldsymbol{\sigma}^2 (I- P\boldsymbol{\phi}_1^2)^{-1}\mathbf{1}
          - (\boldsymbol{\mu} \boldsymbol{\phi}_0 (I- P\boldsymbol{\phi}_1)^{-1}\mathbf{1})^2. \NN
\EEA
This result can be considered as an appropriate upper bound for the variance as we utilize inequality (\ref{ineq E (sum prod a_i)^2}) for the first term of (\ref{E Y2}) by Cauchy Schwarz inequality.
\end{Rem}

\begin{Theo} \label{Theo exist}
    Let $\{Y_t\}_{t=0}^\infty$ follows the HM-MAR($K,1$) model defined by (\ref{HM MAR1}) and $\boldsymbol{\phi}_i^+ = diag(|a_{i,1}|,\, \cdots ,\, |a_{i,K}|)$ for $i=0,1$. \textbf{Also,} let conditions of theorem \ref{Theo second HM-MAR1} hold and all eigenvalues  of $P\boldsymbol{\phi}_1^+$ lie inside the unit circle\textbf{,}
    then $E[\lim{t \rightarrow \infty}Y^2_t]$ exists and is finite.
\end{Theo}

\begin{proof}
    Let random variable $X$ be defined as
    \BEA
      \NewAND  X = \lim_{t \rightarrow \infty}X_t = \lim_{t \rightarrow \infty} \{ \sum_{m=0}^{t} |\prod_{i=2}^{m+1} a_{1,Z_{i}} |(|a_{0,Z_{1}}| +| \sigma_{Z_{1}} \varepsilon_{t-m}| )
             +   |\prod_{i=0}^{t} a_{1,Z_{i}} Y_{0}|  \} \NN.
    \EEA
    By monotone convergence theorem (theorem 16.2 of \cite{Billingsley}) $E[X^2] = \lim_{t \rightarrow \infty}E[X_t^2]$.
    By the assumption of theorem \ref{Theo second HM-MAR1}, we deduce that spectral radius of $\mathbf{1}(P(\boldsymbol{\phi}_1^+)^2\mathbf{1})^\prime \boldsymbol{I}$ lies inside the unit circle, so by a similar method as used to obtain (\ref{ineq E (sum prod a_i)^2}) in lemma \ref{Lem E (sum phi1 phi0)2}, we have
    \BEA
         \lim_{t \rightarrow \infty} E[(\sum_{m=0}^{t-1} |\prod_{i=2}^{m+1} a_{1,Z_{i}} \sigma_{Z_{1}} \varepsilon_{t-m}|)^2] <
    2(\frac{1+\boldsymbol{\mu}^\prime \boldsymbol{\sigma}^2 \mathbf{1}}{(1-\lambda^{1/2})})^2. \label{ineq E (sum prod a_i varepsilon)^2}
    \EEA
    So by (\ref{ineq E (sum prod a_i varepsilon)^2}),lemma \ref{Lem E (sum phi1 phi0)2} and Cauchy Schwarz inequality we have that
    \BEA
        \lim_{t \rightarrow \infty} E[(\sum_{m=0}^{t-1} |\prod_{i=2}^{m+1} a_{1,Z_{i}} a_{0,Z_{1}}|)
        && \hspace{-0.25in}(\sum_{m=0}^{t-1} |\prod_{i=2}^{m+1} a_{1,Z_{i}} \sigma_{Z_{1}}\varepsilon_{t-m}|)]
        \NNL && \hspace{-0.25in} \leq
        2\frac{(1+\boldsymbol{\mu}^\prime \boldsymbol{\sigma}^2 \mathbf{1})(1+\boldsymbol{\mu}^\prime \boldsymbol{\phi}_0^2 \mathbf{1})}{(1-\lambda^{1/2})^2}. \label{Ineq E |sigma|}
    \EEA
    By the assumption, all eigenvalues  of $P\boldsymbol{\phi}_1^+$ lie inside the unit circle, so by a similar method as used to obtain (\ref{lim E sum prod a_i}) we have that
    $$\lim_{t \rightarrow \infty} E[\{\sum_{m=0}^{t-1} |\prod_{i=2}^{m+1} a_{1,Z_{i}} a_{0,Z_{1}}\varepsilon_{t-m}\}|] =
    \sqrt{\pi/2}\boldsymbol{\mu}^\prime \boldsymbol{\phi}_0^+ (I- P\boldsymbol{\phi}_1^+)^{-1}\mathrm{1}
    <\infty.
    $$
    Therefor using inequality (\ref{Ineq E |sigma|}) instead of (\ref{E sum prod phi varepsilon}) in the proof of theorem \ref{Theo second HM-MAR1} , relation (\ref{lim E Y2 ineq}) changes to , we get
    \BEA
        \lim_{t \rightarrow \infty}E[X_t^2] < && \hspace{-0.25in} 2(\frac{1+\boldsymbol{\mu}^\prime \boldsymbol{\phi}_0^2 \mathbf{1}}{1-\lambda^{1/2}})^2
          +  \boldsymbol{\mu} \boldsymbol{\sigma}^2 (I- P\boldsymbol{\phi}_1^2)^{-1}\mathbf{1}
          + \NNL &&  2\frac{(1+\boldsymbol{\mu}^\prime \boldsymbol{\sigma}^2 \mathbf{1})(1+\boldsymbol{\mu}^\prime \boldsymbol{\phi}_0^2 \mathbf{1})}{(1-\lambda^{1/2})^2}. \NN
    \EEA
    Thus $X^2$ is integrable, so $X$ is integrable. Also by triangular inequality we have $|Y_t| < X$ for all $t$ and thus for all $\omega \in \mathbb{R}$, $\lim_{t \rightarrow \infty} Y_t = Y$, where
    $$Y = \{\sum_{m=0}^{\infty} \prod_{i=2}^{m+1} a_{1,Z_{i}} (a_{0,Z_{1}} + \sigma_{Z_{1}} \varepsilon_{0} )
             +   \prod_{i=0}^{\infty} a_{1,Z_{i}} Y_{0}  \}. $$
    So by continuous mapping theorem \cite{Billingsley_conv} we have that $Y^2_t \rightarrow Y^2$ almost surely. Finally $|Y_t| < X$ implies that $|Y_t^2| < 1+X^2$, so by the integrability of $X^2$, and  dominated convergence theorem (theorem 16.4 of \cite{Billingsley}) we conclude that $E[\lim_{t \rightarrow \infty}Y_t]$ exists and
    \BEA
       E[Y^2]= E[\lim_{t \rightarrow \infty}Y^2_t]=\lim_{t \rightarrow \infty} E[Y^2_t] < \infty. \NN
    \EEA
\end{proof}

\section{Simulation} \label{Sec Simulation}
The hidden process $\{Z_t\}$ in (\ref{HM MAR p1}-\ref{HM MAR}) is assumed to follow a first order Markov structure, so HM-MAR can be considered as a generalization of MAR model. Clarifying, MAR model can be considered as HM-MAR model with independent hidden process $\{Z_t\}$. However, HM-MAR model is more complex, using the past observations to determine the next coefficients, and demanding a longer calculation to estimate the parameters and dynamically updating the weighting coefficients.

\input{epsf}
\epsfxsize=3in \epsfysize=1.7in
\begin{figure} \label{FigForecast}
\centerline{$\hspace{-.1in}$\epsffile{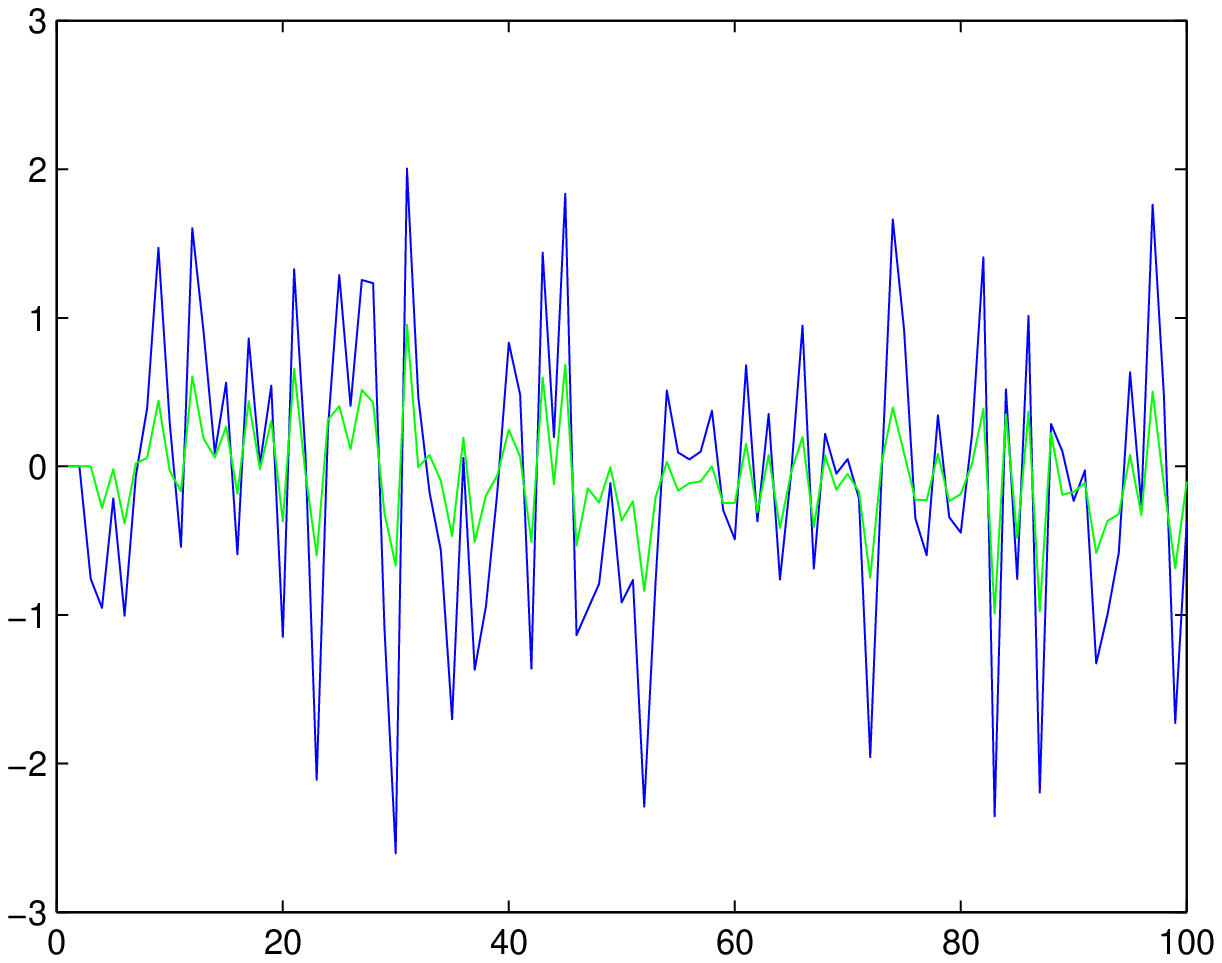}$\hspace{-.22in}$\epsfxsize=3in
\epsfysize=1.7in \epsffile{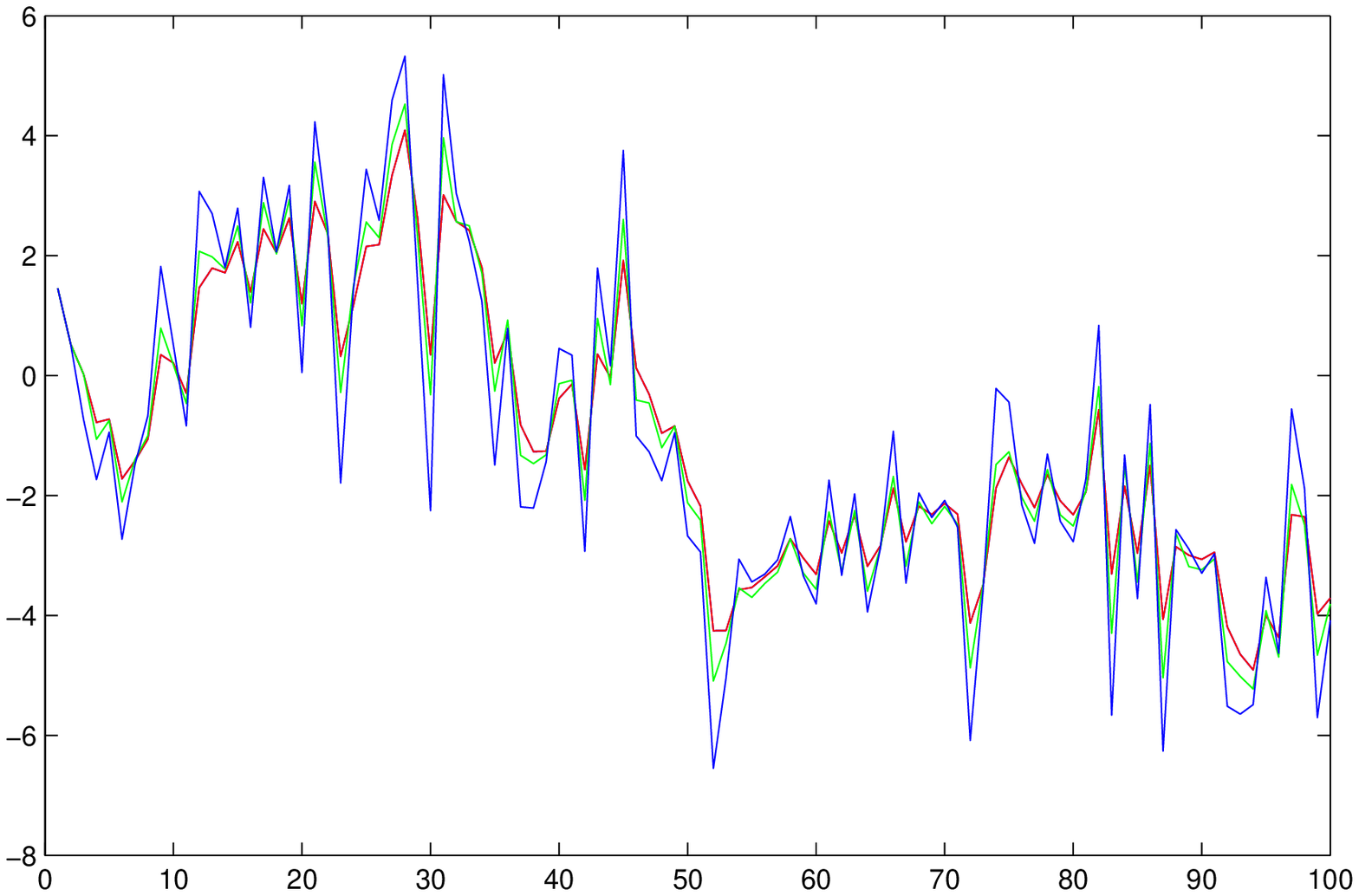}}
 \caption{\scriptsize \hspace{-0.065in} Left: Forecast errors of MAR model (blue), and of HM-MAR model (green). \hspace{4in}$\; \/ \;$ \hspace{0.7in} Right:  Observations of time series (red), forecasts by MAR model (blue), forecasts \hspace{0.1in}$\; \/ \;$ \hspace{0.7in} by HM-MAR model (green).}
\end{figure}

In this section, we investigate the efficiency of these models for time series where $\{Z_t\}$ follow a first order Markov process. To this end, $100$ observations are generated from the following HM-MAR model:
\BEA
    F(y_{t}|\mathcal{F}_{t-1}) & = & \alpha_1^{(t)}
            \Phi(\frac{y_t-0.7y_{t-1}-0.2y_{t-{2}}}{1})
            +
            \NNL &&
            (1-\alpha_1^{(t)})
            \Phi(\frac{y_t-0.5y_{t-1}-0.2y_{t-{2}}}{1}) \nonumber
\EEA
with $\boldsymbol{\rho}=(1,0)'$, (that is starting from the first model ,$\Phi(\frac{y_t-0.7y_{t-1}-0.2y_{t-{2}}}{1})$) , and transition probability matrix $ P= [0.8077 ,\, 0.1923 ;\,  0.7619 ,\, 0.2381]$.
We used EM \cite{mclachlan} algorithm to estimate the conditional probability of hidden variable $Z_t$ given $Y_1,\cdots,Y_T$ (i.e. $P(Z_t|Y_1,\cdots,Y_T)$), and Baum-Welch \cite{macdonald} algorithm to estimate the joint conditional probability of $Z_t,Z_{t-1}$ given $Y_1,\cdots,Y_T$ (i.e. $P(Z_t,Z_{t-1}|Y_1,\cdots,Y_T)$.
Using these estimations we get the following HM-MAR model
\begin{eqnarray*}
    \hat{F}_{HM-MAR}(y_{t}|\mathcal{F}_{t-1}) & = &  \alpha_1^{(t)}
            \Phi(\frac{y_t-0.6514y_{t-1}-0.2973y_{t-{2}}}{0.9887}) +  \\&&
            (1-\alpha_1^{(t)})
            \Phi(\frac{y_t-0.6468y_{t-1}-0.3050y_{t-{2}}}{0.9875}).
\end{eqnarray*}
with  $\boldsymbol{\hat{\rho}}=(0.7261,0.2739)'$ and $ \hat{P}= [ 0.5905, \, 0.4095 ;\,  0.3331 , 0.6669]$, and
\BEA
    \hat{F}_{MAR}(y_{t}|\mathcal{F}_{t-1}) & = &  0.3732
            \Phi(\frac{y_t-0.4042_{t-1}-0.7121y_{t-{2}}}{0.9773}) + \NNL &&
            0.6278
            \Phi(\frac{y_t-0.8176y_{t-1}-0.0485y_{t-{2}}}{0.8640}), \nonumber
\EEA
is the estimated MAR model.
In figure \ref{FigForecast},  the left figure shows the sample path of forecasting errors by MAR model(blue) and forecasting errors by HM-MAR model(green). The right one presents the sample path of simulated HM-MAR model(red), forecasted observations by MAR(blue) and forecasted observations by  HM-MAR model(green). We observe that HM-MAR model produces significantly smaller forecasting errors than MAR model and a better approximation for the time series.
In table \ref{TabIterations} sum of the absolute forecast errors for MAR and HM-MAR models for ten iterations are presented.

\begin{table} \label{TabIterations}
\begin{center}
\caption{Sum of absolute forecasting errors by MAR and HM-MAR models in 10 iterations}\label{tb:margins}
\begin{tabular*}{1\textwidth}%
     {@{\extracolsep{\fill}}ccccccc}
 \hline
   Iterations &&1&2&3&4&5\\
\hline
  HM-MAR &&27.7559&27.7590&27.7560&27.7559&27.7576
 \\
\hline
 MAR &&75.1827&75.0973&75.2079&75.0839&75.2378 \\
\hline
Iterations &&6&7&8&9&10 \\
\hline
HM-MAR && 27.7580&27.7567&27.75657&27.7566&27.7569 \\
\hline
MAR && 75.1054&75.1464&75.1338&75.09641&75.2176
\end{tabular*}
\end{center}
\end{table}

\section{Summary and discussions} \label{Sec Conclusion}

We proposed HM-MAR model as a flexible structure for modeling conditional distribution of $Y_t$ given past observations $(Y_1,\cdots,Y_{t-1})$ in a nonlinear time series. We considered HM-MAR model as the  mixture of some Gaussian distributions \textbf{where} the mean of each distribution follows an AR(p) model.
Unlike the ordinary mixture models, the weighting coefficients \textbf{determining}  the contribution of distributions are not of predefined fixed form (constant values)\textbf{.} \textbf{These values} are conditional probabilities of a latent variable $Z_t$ given past observations $(Y_1,\cdots,Y_{t-1})$.
At each time step $t$, the coefficients are determined through maximizing the posterior probability of latent variable $Z_{t-1}$ given past information.
Latent variables are assumed to follow a Markov process to build a parsimonious model.
 \textbf{A suitable} application \textbf{for} HM-MAR model \textbf{is} when the process $Y_t$ is a result of some processes,
and the contribution of each process changes \textbf{over time}.
  \textbf{If such  effects are} not present in time series \textbf{then} our model automatically will \textbf{reduce} to the ordinary mixture models.
HM-MAR model will also lead to hidden Markov model for continuous process $\{Y_t\}$ \textbf{where} $p$ \textbf{is zero }(i.e. $Y_t$ given $Z_t$, is independent of past observations).

Although modeling the effect of \textbf{all} past information makes the model complicated, a dynamic programming method \textbf{is} proposed for forecasting. \textbf{It is worth mentioning that it} is still possible to study some properties of $\{Y_t\}$, \textbf{such} as
asymptotic behavior of first moment, existence and finiteness of second moment and deriving the upper bound of asymptotic variance of process.
Although the variances of each distribution in the mixture model \textbf{are} constant, the conditional variance of the process in HM-MAR model is not fixed.
This feature can be used to model conditional volatility effects \textbf{ frequently }presented in financial time series.
Another interesting feature is that the first order HM-MAR($K,1$) model can be considered as a mixture of some explosive autoregressive processes (i.e. $a_{.,1} >1$) and the non-explosive ones (i.e. $a_{.,1} < 1$)\textbf{.} \textbf{However, it is} still asymptotically stable in \textbf{first}  and second \textbf{order}.

This \textbf{work} has the potential to be applied in the context of nonlinear time series by imposing hidden Markov property
for the weighting coefficients of mixture model. Also it can elaborate further researches
for extending the stability results to the case of HM-MAR($K,p$), \textbf{where} the lag of autoregressive processes is of order $p$. \textbf{Stationarity }and ergodicity  \textbf{are two major aspects}. Finally \textbf{this area of research can be expanded by considering  other distributions besides} the Gaussian as the underling distribution of mixture model.


\end{document}